\documentclass[a4paper,12pt]{article}
\usepackage{amssymb}
\usepackage[latin1]{inputenc}     

\newcommand{\beq}{\begin{equation} }
\newcommand{\eqq}{\end{equation} }
\newcommand{\cuad}{{\sqcap\kern-.68em\sqcup}}

\newtheorem{remark}{Remark}[section]
\newcommand{\bremark}{\begin{remark} \em}
\newcommand{\eremark}{\end{remark} }

\def\beeq{\begin{equation}}
\def\eeq{\end{equation}}
\newcommand{\begeqaet}{\begin{eqnarray*}}
\newcommand{\eneqaet}{\end{eqnarray*}}

\let\Section=\section
\def\section{\setcounter{equation}{0}\Section}
\newtheorem{Lem}{Lemma}[section]
\newtheorem{Thm}{Theorem}[section]

\newtheorem{Remark}{Remark}[section]

\begin{document}
\begin{center}{\bf\Large Existence of solution for a class of fractional Hamiltonian systems}\medskip

\bigskip

\bigskip

{C\'esar Torres}

 Departamento de Ingenier\'{\i}a  Matem\'atica and
Centro de Modelamiento Matem\'atico
 UMR2071 CNRS-UChile,
 Universidad de Chile\\
 Casilla 170 Correo 3, Santiago, Chile.\\
 {\sl  (ctorres@dim.uchile.cl)}

\end{center}

\medskip

\medskip
\medskip
\medskip
\medskip

\begin{abstract}
In this work we want to prove the existence of solution for a class of fractional Hamiltonian systems given by
\begin{eqnarray*}
_{t}D_{\infty}^{\alpha}(_{-\infty}D_{t}^{\alpha}u(t)) + L(t)u(t)  = & \nabla W(t,u(t))\\
u\in H^{\alpha}(\mathbb{R}, \mathbb{R}^{N})
\end{eqnarray*}
\end{abstract}
\date{}

\setcounter{equation}{0}
\section{ Introduction}

Fractional differential equations both ordinary and partial ones are applied in mathematical modeling of processes in
physics, mechanics, control theory, biochemistry, bioengineering and economics. Therefore the theory of fractional
differential equations is an area intensively developed during last decades \cite{OAJTMJS}, \cite{RH}, \cite{RMJK}, \cite{JSOAJTM}, \cite{BWMBPG}. The monographs \cite{AKHSJT}, \cite{KMBR}, \cite{IP}, enclose a review of methods of solving which are an extension of procedures from differential equations theory.

Recently, also equations including both - left and right fractional derivatives, are discussed. Let us point out that according to integration by parts formulas in fractional calculus , we obtain equations mixing left and right operators. Apart from their possible applications, equations with left and right derivatives are an interesting and new field in fractional differential equations theory. Some works in this topic can be founded in papers \cite{TABS}, \cite{DBJT}, \cite{MK} and its reference.

In this topic recently Jiao and Zhou \cite{FJYZ}, have studied the following fractional boundary value problem
\begin{eqnarray}\label{Eq01}
_{t}D_{T}^{\alpha}({_{0}D_{t}^{\alpha}}u(t)) & = & \nabla F(t,u(t)),\;\mbox{a.e.}\;t\in [0,T],\\
u(0) & = & u(T) = 0.\nonumber
\end{eqnarray}
Using critical point theory they proved the existence of solutions to (\ref{Eq01}).

Motivated by this work, in this paper we want to consider one a class of fractional Hamiltonian systems on $\mathbb{R}$, that is:
\begin{eqnarray}\label{Eq02}
_{t}D_{\infty}^{\alpha}(_{-\infty}D_{t}^{\alpha}u(t)) + L(t)u(t)  = & \nabla W(t,u(t))
\end{eqnarray}
where $\alpha \in (1/2,1)$, $t\in \mathbb{R}$, $u\in \mathbb{R}^{n}$, $L\in C(\mathbb{R}, \mathbb{R}^{n\times n})$ is a symmetric matrix valued function and $W:\mathbb{R}\times \mathbb{R}^{n} \to \mathbb{R}$; satisfies the following condition
\begin{itemize}
\item[$(L)$] $L(t)$ is positive definite symmetric matrix for all $t\in \mathbb{R}$ and there exists an $l\in C(\mathbb{R}, (0,\infty))$ such that $l(t) \to +\infty$ as $t \to \infty$ and
    \begin{equation}\label{Eq03}
    (L(t)x,x) \geq l(t)|x|^{2},\;\;\mbox{for all}\;t\in \mathbb{R}\;\;\mbox{and}\;\;x\in \mathbb{R}^{n}.
    \end{equation}
\item[$(W_{1})$] $W\in C^{1}(\mathbb{R} \times \mathbb{R}^{n}, \mathbb{R})$ and there is a constant $\mu >2$ such that
$$
0< \mu W(t,x) \leq (x, \nabla W(t,x)),\;\;\mbox{for all}\;t\in \mathbb{R}\;\;\mbox{and}\;x\in\mathbb{R}^{n}\setminus \{0\}.
$$
\item[$(W_{2})$] $|\nabla W(t,x)| = o(|x|)$ as $x\to 0$ uniformly with respect to $t\in \mathbb{R}$.
\item[$(W_{3})$] There exists $\overline{W} \in C(\mathbb{R}^{n}, \mathbb{R})$ such that
$$
|W(t,x)| + |\nabla W(t,x)| \leq |\overline{W(x)}|\;\;\mbox{for every}\;x\in \mathbb{R}^{n}\;\mbox{and}\;t\in \mathbb{R}.
$$
\end{itemize}

In particular, if $\alpha = 1$, (\ref{Eq02}) reduces to the standard second order Hamiltonian system
\begin{equation}\label{HEq01}
u'' - L(t)u + \nabla W(t,u)=0,
\end{equation}
where $W: \mathbb{R} \times \mathbb{R}^{n} \to \mathbb{R}$ is a given function and $\nabla W(t,u)$ is the gradient of $W$ at $u$. The existence of homoclinic solution is one of the most important problems in the history of Hamiltonian systems, and has been studied intensively by many mathematicians. Assuming that $L(t)$ and $W(t,u)$ are independent of $t$, or $T$-periodic in $t$, many authors have studied the existence of homoclinic solutions for Hamiltonian systems (\ref{HEq01}) via critical point theory and variational methods. In this case, the existence of homoclinic solution can be obtained by going to the limit of periodic solutions of approximating problems.

If $L(t)$ and $W(t,u)$ are neither autonomous nor periodic in $t$, this problem is quite different from the ones just described, because the lack of compacteness of the Sobolev embedding. In \cite{PRKT} the authors considered (\ref{HEq01}) without periodicity assumptions on $L$ and $W$ and showed that (\ref{HEq01}) possesses one homoclinic solution by using a variant of the mountain pass theorem without the Palais-Smale contidion. In \cite{WOMW}, under the same assumptions of \cite{PRKT}, the authors, by employing a new compact embedding theorem, obtained the existence of homoclinic solution of (\ref{HEq01}).

Our goal in this paper is to show how variational methods based on Mountain pass theorem can be used to get existence results for (\ref{Eq02}). However, the direct application of the mountain pass theorem is not enough since the Palais-Smale sequences might lose compactness in the whole space $\mathbb{R}$. To overcome this difficulty we proof a version of compact embedding for fractional space following the ideas of \cite{WOMW}. Before stating our results let us introduce the main ingredients involved in our approach. We define
$$
\|u\|_{I_{-\infty}^{\alpha}}^{2} = \int_{-\infty}^{\infty} |u(t)|^{2}dt + \int_{-\infty}^{\infty}|_{-\infty}D_{t}^{\alpha}u(t)|^{2}dt
$$
and the space
$$
I_{-\infty}^{\alpha}(\mathbb{R}) = \overline{C_{0}^{\infty}(\mathbb{R}, \mathbb{R}^{n})}^{\|.\|_{\alpha}}.
$$
Now we say that $u\in I_{-\infty}^{\alpha}(\mathbb{R})$ is a weak solution of (\ref{Eq02}) if
$$
\int_{-\infty}^{\infty} [(_{-\infty}D_{t}^{\alpha}u(t), {_{-\infty}}D_{t}^{\alpha}v(t)) + (L(t)u(t),v(t))]dt = \int_{-\infty}^{\infty} (\nabla W(t,u(t)), v(t))dt,
$$
for all $v\in I_{-\infty}^{\alpha}(\mathbb{R})$. For $u\in I_{-\infty}^{\alpha}(\mathbb{R})$ we may define the functional
\begin{equation}\label{Eq03}
I(u) = \frac{1}{2}\int_{-\infty}^{\infty}[|_{-\infty}D_{t}^{\alpha}u(t)|^{2} + (L(t)u(t),u(t))]dt - \int_{-\infty}^{\infty}W(t,u(t))dt.
\end{equation}
which is of class $C^1$. We say that $u\in E^{\alpha}$ is a weak solution of (\ref{Eq02}) if $u$ is a critical point of $I$.

Now we are in a position to state our main existence theorem
\begin{Thm}\label{tm01}
Suppose that $(L), (W_{1})-(W_{2})$ hold, then (\ref{Eq02}) possesses at least one nontrivial solution.
\end{Thm}

The rest of the paper is organized as follows: in section 2, subsection 2.1, we describe the Liouville-Weyl fractional calculus; in subsection 2.2 we introduce the fractional space that we use in our work and some proposition are proven which will aid in our analysis. In section 3, we will prove theorem \ref{tm01}.

\section{Preliminary Results}

\subsection{Liouville-Weyl Fractional Calculus}

The Liouville-Weyl fractional integrals of order $0<\alpha < 1$ are defined as
\begin{equation}\label{LWeq01}
_{-\infty}I_{x}^{\alpha}u(x) = \frac{1}{\Gamma (\alpha)} \int_{-\infty}^{x}(x-\xi)^{\alpha - 1}u(\xi)d\xi
\end{equation}
\begin{equation}\label{LWeq02}
_{x}I_{\infty}^{\alpha}u(x) = \frac{1}{\Gamma (\alpha)} \int_{x}^{\infty}(\xi - x)^{\alpha - 1}u(\xi)d\xi
\end{equation}
The Liouville-Weyl fractional derivative of order $0<\alpha <1$ are defined as the left-inverse operators of the corresponding Liouville-Weyl fractional integrals
\begin{equation}\label{LWeq03}
_{-\infty}D_{x}^{\alpha}u(x) = \frac{d }{d x} {_{-\infty}}I_{x}^{1-\alpha}u(x)
\end{equation}
\begin{equation}\label{LWeq04}
_{x}D_{\infty}^{\alpha}u(x) = -\frac{d }{d x} {_{x}}I_{\infty}^{1-\alpha}u(x)
\end{equation}
The definitions (\ref{LWeq03}) and (\ref{LWeq04}) may be written in an alternative form:
\begin{equation}\label{LWeq05}
_{-\infty}D_{x}^{\alpha}u(x) = \frac{\alpha}{\Gamma (1-\alpha)} \int_{0}^{\infty}\frac{u(x) - u(x-\xi)}{\xi^{\alpha + 1}}d\xi
\end{equation}
\begin{equation}\label{LWeq05}
_{x}D_{\infty}^{\alpha}u(x) = \frac{\alpha}{\Gamma (1-\alpha)} \int_{0}^{\infty}\frac{u(x) - u(x+\xi)}{\xi^{\alpha + 1}}d\xi
\end{equation}

\noindent
We establish the Fourier transform properties of the fractional integral and fractional differential operators. Recall that the Fourier transform $\widehat{u}(w)$ of $u(x)$ is defined by
$$
\widehat{u}(w) = \int_{-\infty}^{\infty} e^{-ix.w}u(x)dx.
$$
Let $u(x)$ be defined on $(-\infty, \infty)$. Then the Fourier transform of the Liouville-Weyl integral and differential operator satisfies
\begin{equation}\label{LWeq06}
\widehat{ _{-\infty}I_{x}^{\alpha}u(x)}(w) = (iw)^{-\alpha}\widehat{u}(w)
\end{equation}
\begin{equation}\label{LWeq07}
\widehat{ _{x}I_{\infty}^{\alpha}u(x)}(w) = (-iw)^{-\alpha}\widehat{u}(w)
\end{equation}
\begin{equation}\label{LWeq08}
\widehat{ _{-\infty}D_{x}^{\alpha}u(x)}(w) = (iw)^{\alpha}\widehat{u}(w)
\end{equation}
\begin{equation}\label{LWeq09}
\widehat{ _{x}D_{\infty}^{\alpha}u(x)}(w) = (-iw)^{\alpha}\widehat{u}(w)
\end{equation}
\subsection{Fractional Derivative Spaces}

In this section we introduce some fractional spaces for more detail see \cite{VEJR}.

\noindent
Let $\alpha > 0$. Define the semi-norm
$$
|u|_{I_{-\infty}^{\alpha}} = \|_{-\infty}D_{x}^{\alpha}u\|_{L^{2}}
$$
and norm
\begin{equation}\label{FDEeq01}
\|u\|_{I_{-\infty}^{\alpha}} = \left( \|u\|_{L^{2}}^{2} + |u|_{I_{-\infty}^{\alpha}}^{2} \right)^{1/2},
\end{equation}
and let
$$
I_{-\infty}^{\alpha} (\mathbb{R}) = \overline{C_{0}^{\infty}(\mathbb{R})}^{\|.\|_{I_{-\infty}^{\alpha}}}.
$$
Now we define the fractional Sobolev space $H^{\alpha}(\mathbb{R})$ in terms of the fourier transform. Let $0< \alpha < 1$, let the semi-norm
\begin{equation}\label{FDEeq02}
|u|_{\alpha} = \||w|^{\alpha}\widehat{u}\|_{L^{2}}
\end{equation}
and norm
$$
\|u\|_{\alpha} = \left( \|u\|_{L^{2}}^{2} + |u|_{\alpha}^{2} \right)^{1/2},
$$
and let
$$
H^{\alpha}(\mathbb{R}) = \overline{C_{0}^{\infty}(\mathbb{R})}^{\|.\|_{\alpha}}.
$$

\noindent
We note a function $u\in L^{2}(\mathbb{R})$ belong to $I_{-\infty}^{\alpha}(\mathbb{R})$ if and only if
\begin{equation}\label{FDEeq03}
|w|^{\alpha}\widehat{u} \in L^{2}(\mathbb{R}).
\end{equation}
Especially
\begin{equation}\label{FDEeq04}
|u|_{I_{-\infty}^{\alpha}} = \||w|^{\alpha}\widehat{u}\|_{L^{2}}.
\end{equation}
Therefore $I_{-\infty}^{\alpha}(\mathbb{R})$ and $H^{\alpha}(\mathbb{R})$ are equivalent with equivalent semi-norm and norm. Analogous to $I_{-\infty}^{\alpha}(\mathbb{R})$ we introduce $I_{\infty}^{\alpha}(\mathbb{R})$. Let the semi-norm
$$
|u|_{I_{\infty}^{\alpha}} = \|_{x}D_{\infty}^{\alpha}u\|_{L^{2}}
$$
and norm
\begin{equation}\label{FDEeq05}
\|u\|_{I_{\infty}^{\alpha}} = \left( \|u\|_{L^{2}}^{2} + |u|_{I_{\infty}^{\alpha}}^{2} \right)^{1/2},
\end{equation}
and let
$$
I_{\infty}^{\alpha}(\mathbb{R}) = \overline{C_{0}^{\infty}(\mathbb{R})}^{\|.\|_{I_{\infty}^{\alpha}}}.
$$
Moreover $I_{-\infty}^{\alpha}(\mathbb{R})$ and $I_{\infty}^{\alpha}(\mathbb{R})$ are equivalent , with equivalent semi-norm and norm \cite{VEJR}.

Now we give the prove of the Sobolev lemma.
\begin{Thm}\label{FDEtm01}
If $\alpha > \frac{1}{2}$, then $H^{\alpha}(\mathbb{R}) \subset C(\mathbb{R})$ and there is a constant $C=C_{\alpha}$ such that
\begin{equation}\label{FDEeq06}
\sup_{x\in \mathbb{R}} |u(x)| \leq C \|u\|_{\alpha}
\end{equation}
\end{Thm}

\noindent
{\bf Proof.} By the Fourier inversion theorem, if $\widehat{u} \in L^{1}(\mathbb{R})$, then $u$ is continuous and
$$
\sup_{x\in \mathbb{R}} |u(x)| \leq \|\widehat{u}\|_{L^{1}}.
$$
Hence, to prove the theorem it is enough to prove that
$$
\|\widehat{u}\|_{L^{1}} \leq \|u\|_{\alpha},
$$
so by Schwarz inequality, we have
\begin{eqnarray*}
\int_{\mathbb{R}} |\widehat{u}(w)|dw & = & \int_{\mathbb{R}} (1 + |w|^{2})^{\alpha /2}|\widehat{u}(w)| \frac{1}{(1+|w|^{2})^{\alpha /2}}dw\\
                                     & \leq & \left( \int_{\mathbb{R}} (1 + |w|^{2\alpha}) |\widehat{u}(w)|^{2}dw \right)^{1/2}\left( \int_{\mathbb{R}}(1+|w|^{2})^{-\alpha}dw \right)^{1/2}.
\end{eqnarray*}
The first integral on the right is $\|u\|_{\alpha}^{2}$, so the theorem boils down to the fact
$$
\int_{\mathbb{R}} (1 + |w|^{2})^{-\alpha}dw = \int_{0}^{\infty} (1+r^{2})^{-\alpha}r^{n-1}dr < \infty
$$
precisely when $\alpha > \frac{1}{2}$. $\Box$
\begin{Remark}\label{FDEnta01}
If $u\in H^{\alpha}(\mathbb{R})$, then $u\in L^{q}(\mathbb{R})$ for all $q\in [2,\infty]$, since
$$
\int_{\mathbb{R}} |u(x)|^{q}dx \leq \|u\|_{\infty}^{q-2}\|u\|_{L^{2}}^{2}
$$
\end{Remark}

\noindent
Now we introduce a new fractional spaces. Let
$$
X^{\alpha} = \left\{ u\in H^{\alpha}(\mathbb{R}, \mathbb{R}^{n})|\;\;\int_{\mathbb{R}} |_{-\infty}D_{t}^{\alpha}u(t)|^{2} + L(t)u(t).u(t) dt < \infty  \right\}
$$
The space $X^{\alpha}$ is a Hilbert space with the inner product
$$
\langle u,v \rangle_{X^{\alpha}} = \int_{\mathbb{R}} (_{-\infty}D_{t}^{\alpha}u(t) , \; _{-\infty}D_{t}^{\alpha}v(t)) + L(t)u(t).v(t)dt
$$
and the corresponding norm
$$
\|u\|_{X^{\alpha}}^{2} = \langle u,u \rangle_{X^{\alpha}}
$$
\begin{Lem}\label{FDElm01}
Suppose $L$ satisfies ($L$). Then $X^{\alpha}$ is continuously embedded in $H^{\alpha}(\mathbb{R},\mathbb{R}^{n})$.
\end{Lem}

\noindent
{\bf Proof.} Since $l\in C(\mathbb{R}, (0,\infty))$ and $l$ is coercive, then $l_{min} = \min_{t\in \mathbb{R}}l(t)$ exists, so we have
$$
(L(t)u(t) , u(t)) \geq l(t)|u(t)|^{2} \geq l_{min}|u(t)|^{2},\;\forall t\in \mathbb{R}.
$$
Then
\begin{eqnarray*}
l_{min}\|u\|_{\alpha}^{2} & = & l_{min}\left( \int_{\mathbb{R}} |_{-\infty}D_{t}^{\alpha}u(t)|^{2} + |u(t)|^{2}dt\right)\\
                          & \leq & l_{min} \int_{\mathbb{R}}|_{-\infty}D_{t}^{\alpha}u(t)|^{2}dt + \int_{\mathbb{R}}(L(t)u(t),u(t))dt
\end{eqnarray*}
So
\begin{equation}\label{FDEeq07}
\|u\|_{\alpha}^{2} \leq K \|u\|_{X^{\alpha}}^{2}
\end{equation}
where $K = \frac{\max \{l_{min}, 1\}}{l_{min}}$. $\Box$
\begin{Lem}\label{FDElm02}
Suppose $L$ satisfies ($L$). Then the imbedding of $X^{\alpha}$ in $L^{2}(\mathbb{R})$ is compact.
\end{Lem}

\noindent
{\bf Proof.} We note first that by lemma \ref{FDElm01} and remark \ref{FDEnta01} we have
$$
X^{\alpha} \hookrightarrow L^{2}(\mathbb{R})\;\;\mbox{is continuous}.
$$
Now, let $(u_{k}) \in X^{\alpha}$ be a sequence such that $u_{k} \rightharpoonup u$ in $X^{\alpha}$. We will show that $u_{k} \to u$ in $L^{2}(\mathbb{R})$. Suppose, without loss of generality, that $u_{k} \to 0$ in $X^{\alpha}$. The Banach-Steinhaus theorem implies that
$$
A = \sup_{k}\|u_{k}\|_{X^{\alpha}} < +\infty
$$
Let $\epsilon >0$; there is $T_{0}<0$ such that $\frac{1}{l(t)} \leq \epsilon$ for all $t$ such that $t\leq T_{0}$. Similarly, there is $T_{1}>0$, such that $\frac{1}{l(t)}\leq \epsilon$ for all $t\geq T_{1}$. Sobolev's theorem (see e.g. \cite{CAS}) implies that $u_{k} \to 0$ uniformly on $\overline{\Omega} = [T_{0}, T_{1}]$, so there is a $k_{0}$ such that
\begin{equation}\label{FDEeq08}
\int_{\Omega} |u_{k}(t)|^{2}dt \leq \epsilon,\;\;\mbox{for all}\;k\geq k_{0}.
\end{equation}
Since $\frac{1}{l(t)} \leq \epsilon$ on $(-\infty , T_{0}]$ we have
\begin{equation}\label{FDEeq09}
\int_{-\infty}^{T_{0}} |u_{k}(t)|^{2}dt \leq \epsilon \int_{-\infty}^{T_{0}} l(t)|u_{k}(t)|^{2}dt \leq \epsilon A^{2}.
\end{equation}
Similarly, since $\frac{1}{l(t)} \leq \epsilon$ on $[T_{1}, +\infty)$, we have
\begin{equation}\label{FDEeq10}
\int_{T_{1}}^{+\infty} |u_{k}(t)|^{2}dt \leq \epsilon A^{2}.
\end{equation}
Combining (\ref{FDEeq08}), (\ref{FDEeq09}) and (\ref{FDEeq10}) we get $u_{k} \to 0$ in $L^{2}(\mathbb{R}, \mathbb{R}^{n})$. $\Box$

\begin{Lem}\label{FDE-lem01}
There are constants $c_{1}>0$ and $c_{2}>0$ such that
\begin{equation}\label{FDE-eq01}
W(t,u) \geq c_{1}|u|^{\mu},\;\;|u|\geq 1
\end{equation}
and
\begin{equation}\label{FDE-eq02}
W(t,u) \leq c_{2}|u|^{\mu},\;\;|u|\leq 1
\end{equation}
\end{Lem}

{\bf Proof.} By $(W_{1})$ we note that
$$
\mu W(t,\sigma u) \leq (\sigma u, \nabla W(t, \sigma u)).
$$
Let $f(\sigma) = W(t,\sigma u)$, then
\begin{equation}\label{FDE-eq03}
\frac{d}{d\sigma}\left( f(\sigma) \sigma^{-\mu} \right)\geq 0
\end{equation}
Now we consider two cases

\noindent
{\bf Case 1.} $|u|\leq 1$. In this case we integrate (\ref{FDE-eq03}), from $1$ until $\frac{1}{|u|}$ and we get
\begin{equation}\label{FDE-eq04}
W(t,u) \leq W(t,\frac{u}{|u|})|u|^{\mu}.
\end{equation}

\noindent
{\bf Case 2.} $|u|\geq 1$. In this case we integrate (\ref{FDE-eq03}), from $\frac{1}{|u|}$ until $1$ and we get
\begin{equation}\label{FDE-eq05}
W(t,u) \geq |u|^{\mu}W(t, \frac{u}{|u|}).
\end{equation}
Now, since $u\in \mathbb{R}^{n}$, $\frac{u}{|u|}\in B(0,1)$. So, since $W$ is continuous and $B(0,1)$ is compact, there are $c_{1}>0$ and $c_{2}>0$ such that
$$
c_{1} \leq W(t,u) \leq c_{2},\;\;\mbox{for every}\;u\in B(0,1).
$$
Therefore we get the affirmation of the lemma. $\Box$

\begin{Remark}\label{FDE-nta01}

By lemma \ref{FDE-lem01}, we have
\begin{equation}\label{FDE-eq06}
W(t, u) = o(|u|^{2})\;\mbox{as}\;u\to 0\;\mbox{uniformly in}\;t\in\mathbb{R}
\end{equation}
In addition, by $(W_{2})$, we have, for any $u\in \mathbb{R}^{n}$ such that $|u| \leq M_{1}$, there exists some constant $d>0$ (dependent on $M_{1}$) such that
\begin{equation}\label{FDE-eq07}
|\nabla W(t,u(t))| \leq d|u(t)|
\end{equation}
\end{Remark}

\noindent
Similar to lemma 2 of \cite{WOMW}, we can get the following result.
\begin{Lem}\label{FDElm03}
Suppose that (L), ($W_{1}$)-($W_{2}$) are satisfied. If $u_{k} \rightharpoonup u$ in $X^{\alpha}$, then $\nabla W(t, u_{k}) \to \nabla W(t, u)$ in $L^{2}(\mathbb{R}, \mathbb{R}^{n})$.
\end{Lem}

\noindent
{\bf Proof.} Assume that $u_{k} \rightharpoonup u$ in $X^{\alpha}$. Then there exists a constant $d_{1}>0$ such that, by Banach-Steinhaus theorem and (\ref{FDEeq06}),
$$
\sup_{k\in \mathbb{N}} \|u_{k}\|_{\infty} \leq d_{1},\;\;\|u\|_{\infty} \leq d_{1}.
$$
By ($W_{2}$), for any $\epsilon >0$ there is $\delta >0$ such that
$$
|u_{k}|< \delta\;\;\mbox{implies}\;\;|\nabla W(t,u_{k})|\leq \epsilon |u_{k}|
$$
and by ($W_{3}$) there is $M>0$ such that
$$
|\nabla W(t,u_{k})| \leq M,\;\;\mbox{forall}\;\delta < u_{k} \leq d_{1}.
$$
Therefore, there exists a constant $d_{2}>0$ (dependent on $d_{1}$) such that
$$
|\nabla W(t, u_{k}(t))| \leq d_{2}|u_{k}(t)|,\;\;|\nabla W(t,u(t))| \leq d_{2}|u(t)|
$$
for all $k\in \mathbb{N}$ and $t\in \mathbb{R}$. Hence,
$$
|\nabla W(t,u_{k}(t)) - \nabla W(t,u(t))| \leq d_{2}(|u_{k}(t)| + |u(t)|) \leq d_{2}(|u_{k}(t) - u(t)| + 2|u(t)|),
$$
Since, by lemma \ref{FDElm02}, $u_{k} \to u$ in $L^{2}(\mathbb{R}, \mathbb{R}^{n})$, passing to a subsequence if necessary, it can be assumed that
$$
\sum_{k=1}^{\infty} \|u_{k} - u\|_{L^{2}} < \infty
$$
But this implies $u_{k}(t) \to u(t)$ almost everywhere $t\in \mathbb{R}$ and
$$
\sum_{k=1}^{\infty}|u_{k}(t) - u(t)| = v(t) \in L^{2}(\mathbb{R},\mathbb{R}^{n}).
$$
Therefore
$$
|\nabla W(t,u_{k}(t)) - \nabla W(t,u(t))| \leq d_{2}(v(t) + 2|u(t)|).
$$
Then, using the Lebesgue's convergence theorem, the lemma is proved. $\Box$

Now we introduce more notations and some necessary definitions. Let $\mathfrak{B}$ be a real Banach space, $I\in C^{1}(\mathfrak{B},\mathbb{R})$, which means that $I$ is a continuously Fréchet-differentiable functional defined on $\mathfrak{B}$. Recall that $I\in C^{1}(\mathfrak{B},\mathbb{R})$ is said to satisfy the (PS) condition if any sequence $\{u_{k}\}_{k\in \mathbb{N}} \in \mathfrak{B}$, for which $\{I(u_{k})\}_{k\in \mathbb{N}}$ is bounded and $I'(u_{k}) \to 0$ as $k\to +\infty$, possesses a convergent subsequence in $\mathfrak{B}$.

Moreover, let $B_{r}$ be the open ball in $\mathfrak{B}$ with the radius $r$ and centered at $0$ and $\partial B_{r}$ denote its boundary. We obtain the existence of homoclinic solutions of (\ref{Eq02}) by use of the following well-known Mountain Pass Theorems, see \cite{PR}.
\begin{Thm}\label{FDEtm02}
Let $\mathfrak{B}$ be a real Banach space and $I\in C^{1}(\mathfrak{B}, \mathbb{R})$ satisfying (PS) condition. Suppose that $I(0) = 0$ and
\begin{itemize}
\item[i.] There are constants $\rho , \beta >0$ such that $I|_{\partial B_{\rho}} \geq \beta$, and
\item[ii.] There is and $e\in \mathfrak{B} \setminus \overline{B_{\rho}}$ such that $I(e)\leq 0$.
\end{itemize}
Then $I$ possesses a critical value $c\geq \beta$. Moreover $c$ can be characterized as
$$
c = \inf_{\gamma \in \Gamma} \max_{s\in [0,1]} I(\gamma (s))
$$
where
$$
\Gamma = \{\gamma \in C([0,1] , \mathfrak{B}):\;\;\gamma (0) = 0,\;\;\gamma (1) = e\}
$$
\end{Thm}
\section{Proof of Theorem \ref{tm01}}

Now we are going to establish the corresponding variational framework to obtain the existence of solutions for (\ref{Eq02}). Define the functional $I: X^{\alpha} \to \mathbb{R}$ by
\begin{eqnarray}\label{TMeq01}
I(u) & = & \int_{\mathbb{R}} \left[ \frac{1}{2}|_{-\infty}D_{t}^{\alpha}u(t)|^{2} + \frac{1}{2}(L(t)u(t),u(t)) - W(t,u(t))\right]dt \nonumber\\
     & = & \frac{1}{2}\|u\|_{X^{\alpha}}^{2} - \int_{\mathbb{R}} W(t,u(t))dt
\end{eqnarray}
\begin{Lem}\label{TMlm01}
Under the conditions of Theorem \ref{tm01}, we have
\begin{equation}\label{TMeq02}
I'(u)v = \int_{\mathbb{R}} \left[ (_{-\infty}D_{t}^{\alpha}u(t), _{-\infty}D_{t}^{\alpha}v(t)) + (L(t)u(t),v(t)) - (\nabla W(t,u(t)),v(t)) \right]dt
\end{equation}
for all $u,v \in X^{\alpha}$, which yields that
\begin{equation}\label{TMeq03}
I'(u)u = \|u\|_{X^{\alpha}}^{2} - \int_{\mathbb{R}}(\nabla W(t,u(t)), u(t))dt.
\end{equation}
Moreover, $I$ is a continuously Fréchet-differentiable functional defined on $X^{\alpha}$, i.e., $I\in C^{1}(X^{\alpha}, \mathbb{R})$.
\end{Lem}

\noindent
{\bf Proof.} We firstly show that $I: X^{\alpha} \to \mathbb{R}$. By (\ref{FDE-eq06}), there is a $\delta >0$ such that $|u| \leq \delta$ implies that
\begin{equation}\label{TMeq04}
W(t,u) \leq \epsilon |u|^{2}\;\;\mbox{for all}\;t\in \mathbb{R}
\end{equation}
Let $u\in X^{\alpha}$, then $u \in C(\mathbb{R},\mathbb{R}^{n})$, the space of continuous function $u\in \mathbb{R}$ such that $u(t)\to 0$ as $|t| \to +\infty$. Therefore there is a constant $R>0$ such that $|t| \geq R$ implies $|u(t)|\leq \delta$. Hence, by (\ref{TMeq04}), we have
\begin{equation}\label{TMeq05}
\int_{\mathbb{R}} W(t, u(t)) \leq \int_{-R}^{R} W(t,u(t))dt + \epsilon \int_{|t|\geq R}|u(t)|^{2}dt < +\infty.
\end{equation}
Combining (\ref{TMeq01}) and (\ref{TMeq05}), we show that $I:X^{\alpha} \to \mathbb{R}$.

\noindent
Now we prove that $I\in C^{1}(X^{\alpha}, \mathbb{R})$. Rewrite $I$ as follows
$$
I = I_{1} - I_{2},
$$
where
$$
I_{1} = \frac{1}{2} \int_{\mathbb{R}} [|_{-\infty}D_{t}^{\alpha}u(t)|^{2} + (L(t)u(t),u(t))]dt,\;\;I_{2} = \int_{\mathbb{R}}W(t,u(t))dt
$$
It is easy to check that $I_{1} \in C^{1}(X^{\alpha},\mathbb{R})$ and
\begin{equation}\label{TMeq06}
I'_{1}(u)v = \int_{\mathbb{R}}\left[ (_{-\infty}D_{t}^{\alpha}u(t), _{-\infty}D_{t}^{\alpha}v(t)) + (L(t)u(t), v(t))\right]dt.
\end{equation}
Thus it is sufficient to show this is the case for $I_{2}$. In the process we will see that
\begin{equation}\label{TMeq07}
I'_{2}(u)v = \int_{\mathbb{R}}(\nabla W(t,u(t)), v(t))dt,
\end{equation}
which is defined for all $u,v\in X^{\alpha}$. For any given $u\in X^{\alpha}$, let us define $J(u): X^{\alpha} \to \mathbb{R}$ as follows
$$
J(u)v = \int_{\mathbb{R}} (\nabla W(t,u(t)), v(t))dt,\;\;\forall v\in X^{\alpha}.
$$
It is obvious that $J(u)$ is linear. Now we show that $J(u)$ is bounded. Indeed, for any given $u\in X^{\alpha}$, by (\ref{FDE-eq07}), there is a constant $d_{3}>0$ such that
$$
|\nabla W(t,u(t))| \leq d_{3}|u(t)|,
$$
which yields that, by the Hölder inequality and lemma \ref{FDElm01}
\begin{eqnarray}\label{TMeq08}
|J(u)v| & = & \left| \int_{\mathbb{R}} (\nabla W(t, u(t)), v(t)) dt\right| \leq d_{3}\int_{\mathbb{R}}|u(t)||v(t)|dt \nonumber\\
        & \leq & \frac{d_{3}}{l_{min}}\|u\|_{X^{\alpha}}\|v\|_{X^{\alpha}}.
\end{eqnarray}
Moreover, for $u$ and $v \in X^{\alpha}$, by Mean Value theorem, we have
$$
\int_{\mathbb{R}} W(t,u(t) + v(t))dt - \int_{\mathbb{R}} W(t, u(t))dt = \int_{\mathbb{R}} (\nabla W(t,u(t) + h(t)v(t)))dt,
$$
where $h(t)\in (0,1)$. Therefore, by lemma \ref{FDElm02} and the Hölder inequality, we have
\begin{eqnarray}\label{TMeq09}
\int_{\mathbb{R}} (\nabla W(t,u(t) + h(t)v(t)), v(t))dt - \int_{\mathbb{R}} (\nabla W(t, u(t)), v(t))dt \nonumber\\
= \int_{\mathbb{R}} (\nabla W(t,u(t)) + h(t)v(t) - \nabla W(t,u(t)), v(t))dt \to 0
\end{eqnarray}
as $v\to 0$ in $X^{\alpha}$. Combining (\ref{TMeq08}) and (\ref{TMeq09}), we see that (\ref{TMeq07}) holds. It remains to prove that $I'_{2}$ is continuous. Suppose that $u \to u_{0}$ in $X^{\alpha}$ and note that
\begin{eqnarray*}
\sup_{\|v\|_{X^{\alpha}} = 1} |I'_{2}(u)v - I'_{2}(u_{0})v| & = & \sup_{\|v\|_{X^{\alpha}}= 1} \left| \int_{\mathbb{R}} (\nabla W(t,u(t)) - \nabla W(t,u_{0}(t)), v(t))dt \right|\\
& \leq & \sup_{\|v\|_{X^{\alpha}} = 1} \|\nabla W(.,u(.)) - \nabla W(.,u_{0}(.))\|_{L^{2}}\|v\|_{L^{2}}\\
& \leq & \frac{1}{\sqrt{l_{\min}}} \|\nabla W(.,u(.)) - \nabla W(.,u_{0}(.))\|_{L^{2}}
\end{eqnarray*}
By lemma \ref{FDElm02}, we obtain that $I'_{2}(u)v - I'_{2}(u_{0})v \to 0$ as $\|u\|_{X^{\alpha}} \to \|u_{0}\|_{X^{\alpha}}$ uniformly with respect to $v$, which implies the continuity of $I'_{2}$ and $I\in C^{1}(X^{\alpha}, \mathbb{R})$. $\Box$
\begin{Lem}\label{TMlm02}
Under the conditions of $(L)-(W_{2})$, $I$ satisfies the (PS) condition.
\end{Lem}

\noindent
{\bf Proof.} Assume that $(u_{k})_{k\in \mathbb{N}} \in X^{\alpha}$ is a sequence such that $\{I(u_{k})\}_{k\in \mathbb{N}}$ is bounded and $I'(u_{k}) \to 0$ as $k \to +\infty$. Then there exists a constant $C_{1}>0$ such that
\begin{equation}\label{TMeq10}
|I(u_{k})|\leq C_{1},\;\;\|I'(u_{k})\|_{(X^{\alpha})^{*}} \leq C_{1}
\end{equation}
for every $k\in \mathbb{N}$.

\noindent
We firstly prove that $\{u_{k}\}_{k\in \mathbb{N}}$ is bounded in $X^{\alpha}$. By (\ref{TMeq01}), (\ref{TMeq03}) and ($W_{1}$), we have
\begin{eqnarray}\label{TMeq11}
C_{1} + \|u_{k}\|_{X^{\alpha}} & \geq & I(u_{k})-\frac{1}{\mu}I'(u_{k})u_{k} \nonumber\\
& = &\left( \frac{\mu}{2}  - 1\right)\|u_{k}\|_{X^{\alpha}}^{2} - \int_{\mathbb{R}} [W(t,u_{k}(t)) - \frac{1}{\mu}(\nabla W(t,u_{k}(t)), u_{k}(t))]dt \nonumber\\
&\geq & \left( \frac{\mu}{2}  - 1\right)\|u_{k}\|_{X^{\alpha}}^{2}.
\end{eqnarray}
Since $\mu > 2$, the inequality (\ref{TMeq11}) shows that $\{u_{k}\}_{k\in \mathbb{N}}$ is bounded in $X^{\alpha}$. So passing to a subsequence if necessary, it can be assumed that $u_{k} \rightharpoonup u$ in $X^{\alpha}$ and hence, by lemma \ref{FDElm02}, $u_{k} \to u$ in $L^{2}(\mathbb{R},\mathbb{R}^{n})$. It follows from the definition of $I$ that
\begin{eqnarray}\label{TMeq12}
(I'(u_{k}) - I'(u))(u_{k} - u) & &  \nonumber\\
= \|u_{k} - u\|_{X^{\alpha}}^{2} - \int_{\mathbb{R}}[\nabla W(t,u_{k}) - \nabla W(t,u)](u_{k} - u)dt.
\end{eqnarray}
Since $u_{k} \to u$ in $L^{2}(\mathbb{R}, \mathbb{R}^{n})$, we have (see lemma \ref{FDElm03}) $\nabla W(t, u_{k}(t)) \to \nabla W(t,u(t))$ in $L^{2}(\mathbb{R}, \mathbb{R}^{n})$. Hence
$$
\int_{\mathbb{R}} ( \nabla W(t,u_{k}(t))-\nabla W(t,u(t)), u_{k}(t)-u(t) )dt \to 0
$$
as $k\to +\infty$. So (\ref{TMeq12}) implies
$$
\|u_{k} - u\|_{X^{\alpha}} \to 0\;\;\mbox{as}\;\;k\to +\infty.
$$
$\Box$

\noindent
Now we are in the position to give the proof of theorem \ref{tm01}. We divide the proof into several steps.

\noindent
{\bf Proof of theorem \ref{tm01}.}

\noindent
{\bf Step 1.} It is clear that $I(0) = 0$ and $I\in C^{1}(X^{\alpha}, \mathbb{R})$ satisfies the (PS) condition by lemma \ref{TMlm01} and \ref{TMlm02}.

\noindent
{\bf Step 2.} Now We show that there exist constant $\rho >0$ and $\beta >0$ such that $I$ satisfies the condition (i) of theorem \ref{FDEtm02}. By lemma \ref{FDElm02}, there is a $C_{0}> 0$ such that
$$
\|u\|_{L^{2}} \leq C_{0} \|u\|_{X^{\alpha}}.
$$
On the other hand by theorem \ref{FDEtm01}, there is $C_{\alpha}>0$ such that
$$
\|u\|_{\infty} \leq C_{\alpha} \|u\|_{X^{\alpha}}.
$$
By (\ref{FDE-eq06}), for all $\epsilon >0$, there exists $\delta >0$ such that
$$
W(t,u(t)) \leq \epsilon |u(t)|^{2}\;\;\mbox{wherever}\;\;|u(t)| < \delta.
$$
Let $\rho = \frac{\delta}{C_{\alpha}}$ and $\|u\|_{X^{\alpha}} \leq \rho$; we have $\|u\|_{\infty} \leq \frac{\delta}{C_{\alpha}}.C_{\alpha} = \delta$. Hence
$$
|W(t,u(t))| \leq \epsilon |u(t)|^{2}\;\;\mbox{for all}\;t\in \mathbb{R}.
$$
Integrating on $\mathbb{R}$, we get
$$
\int_{\mathbb{R}}W(t,u(t))dt \leq \epsilon \|u\|_{L^{2}}^{2} \leq \epsilon C_{0}^{2}\|u\|_{X^{\alpha}}^{2}
$$
So, if $\|u\|_{X^{\alpha}} = \rho$, then
$$
I(u) = \frac{1}{2}\|u\|_{X^{\alpha}}^{2} - \int_{\mathbb{R}}W(t,u(t))dt \geq (\frac{1}{2} - \epsilon C_{0}^{2})\|u\|_{X^{\alpha}}^{2} = (\frac{1}{2} - \epsilon C_{0}^{2})\rho^{2}.
$$
And it suffices to choose $\epsilon = \frac{1}{4C_{0}^{2}}$ to get
\begin{equation}\label{TMeq13}
I(u) \geq \frac{\rho^{2}}{4C_{0}^{2}} = \beta >0
\end{equation}

\noindent
{\bf Step 3.} It remains to prove that there exists an $e\in X^{\alpha}$ such that $\|e\|_{X^{\alpha}}> \rho$ and $I(e)\leq 0$, where $\rho$ is defined in Step 2. Consider
$$
I(\sigma u) = \frac{\sigma^{2}}{2}\|u\|_{X^{\alpha}}^{2} - \int_{\mathbb{R}}W(t,\sigma u(t))dt
$$
for all $\sigma \in \mathbb{R}$. By (\ref{FDE-eq01}), there is $c_{1}>0$ such that
\begin{equation}\label{TMeq14}
W(t,u(t)) \geq c_{1}|u(t)|^{\mu}\;\;\mbox{for all}\; |u(t)| \geq 1.
\end{equation}
Take some $u\in X^{\alpha}$ such that $\|u\|_{X^{\alpha}} = 1$. Then there exists a subset $\Omega$ of positive measure of $\mathbb{R}$ such that $u(t) \neq 0$ for $t\in \Omega$. Take $\sigma >0$ such that $\sigma |u(t)| \geq 1$ for $t\in \Omega$. Then by (\ref{TMeq14}), we obtain
\begin{equation}\label{TMeq15}
I(\sigma u) \leq \frac{\sigma^{2}}{2} - c_{1}\sigma^{\mu} \int_{\Omega}|u(t)|^{\mu} dt.
\end{equation}
Since $c_{1}>0$ and $\mu >2$, (\ref{TMeq15}) implies that $I(\sigma u) <0$ for some $\sigma >0$ with $\sigma |u(t)|\geq 1$ for $t\in \Omega$ and $\|\sigma u\|_{X^{\alpha}}> \rho$, where $\rho$ is defined in Step 2. By theorem \ref{FDEtm02}, $I$ possesses a critical value $c\geq \beta >0$ given by
$$
c = \inf_{\gamma \in \Gamma}\max_{s\in [0,1]}I(\gamma (s))
$$
where
$$
\Gamma = \{\gamma \in C([0,1], X^{\alpha}):\;\;\gamma (0) = 0,\;\gamma (1) = e\}.
$$
Hence there is $u\in X^{\alpha}$ such that
$$
I(u) = c,\;\;I'(u) = 0
$$
$\Box$



\noindent {\bf Acknowledgements:}
C.T. was  partially supported by MECESUP 0607.

\end{document}